\documentclass[11pt,reqno]{amsart}
\usepackage{amsmath,amsfonts,amsthm,amssymb}
\usepackage[a4paper]{geometry}
\title[Energy inequalities]{Energy inequalities and dispersive estimates for wave equations with time-dependent coefficients}
\author{Jens Wirth}
\address{Jens Wirth, Department of Mathematics, Imperial College London, 180 Queen's Gate, London SW7 2AZ, UK}
\email{j.wirth@imperial.ac.uk}
\def\R{\mathbb R}
\def\E{\mathbb E}
\def\d{\mathrm d}
\def\D{\mathrm D}
\newtheorem{thm}{Theorem}
\theoremstyle{remark}
\newtheorem{expl}{Example}
\newtheorem{probl}{Open Problem}
\begin{document}
\begin{abstract}
We consider wave models with lower order terms and recollect some recent results on energy and dispersive estimates for their solution based on symbolic type estimates for coefficients and partly stabilisation conditions. The exposition is complemented by a collection of open problems. 
\end{abstract}

\maketitle

\noindent {\bf Model problem.} We consider on $\R_t\times \R^n_x$ the Cauchy problem
\begin{subequations}\label{eq:CP}
\begin{align}
 &u_{tt} - a^2(t) \Delta u+2b(t) u_t  + m^2(t) u = 0,\\
 &u(0,\cdot)=u_1,\quad u_t(0,\cdot)=u_2
 \end{align}
 \end{subequations}
for a wave equation with time-dependent coefficients and ask for the interplay between properties of coefficients (and their long-time behaviour) and large-time asymptotics of solutions (measured in terms of energies or certain norms). There is an extensive literature about results; we will sketch some of them and focus on related open problems. 

Problem \eqref{eq:CP} looks far more general than it is. Using Liouville transforms it is possible to reduce the problem to equations with constant $a(t)$ and vanishing $b(t)$. We will not pursue this further and merely restrict our consideration to equations with variable speed and equations with constant speed and variable dissipation. 

\bigskip
\noindent{\bf I. Some classical results on bounded variable propagation speed.} 

\medskip
\noindent{\bf I.($\alpha$)} For the sake of completeness we start our consideration with the wave equation
\begin{equation}
   u_{tt}-\Delta u = 0.
\end{equation}
Then it is well-known that the corresponding energy
\begin{equation}
\E(u;t) = \frac12 \int \big( |\nabla u(t,x)|^2+|u_t(t,x)|^2\big) \d x
\end{equation}
is preserved, $\E(u;t)=\E(u;0)$. Furthermore, if we replace the $L^2$-norm by an $L^q$-norm, $q>2$, we gain decay. This is known as Strichartz type decay estimate or dispersive estimate and goes back to classical papers of Strichartz, \cite{Strichartz:1970}, Pecher, \cite{Pecher:1976}, and von Wahl, \cite{Wahl:1971}. It plays an essential r\^ole in considerations of nonlinear wave equations. In a precise formulation the result reads as
\begin{equation}
  \| \nabla u(t,\cdot) \|_{L^q} + \| u_t(t,\cdot)\|_{L^q} \le C_{pq} (1+t)^{-\frac{n-1}2 (\frac1p-\frac1q)} \big( \| \langle\D\rangle^{r_p+1} u_1\|_{L^p} + \| \langle\D\rangle^{r_p} u_2\|_{L^p} \big)
\end{equation}
for $p\in(1,2]$, $pq=p+q$ and $r_p = n(1/p-1/q)$. As usual we denote $\langle\D\rangle = (\mathrm I-\Delta)^{1/2}$.

There exist several alternative formulations of this estimate, we want to point out two of them. First, if we replace $(1+t)$
by $t$ (and thus allow the decay function to become singular in $t=0$) the regularity can be sharpened to $r_p=(n+1)(1/p-1/q)/2$. Furthermore, this statement is then equivalent to the Strichartz' estimate 
\begin{equation}
\|u\|_{L^q_tL^r_x} + \|\nabla u\|_{L^\infty_t L^2_x} + \| u_t\|_{L^\infty_t L^2_x} \le C
\big( \|\nabla u_1\|_{L^2}+\|u_2\|_{L^2} 
\big) 
\end{equation}
with $2\le q\le\infty$, $2\le r\le\infty$, $1/q+n/r=n/2$ and $1/q+(n-1)/(2r)\le (n-1)/4$ as pointed out in \cite{Keel:1998}.

\medskip
\noindent{\bf I.($\beta$)} The simplest case of time-dependent propagation speed is a periodic coefficient,
\begin{equation}\label{eq:CPbeta}
  u_{tt} - a^2(t) \Delta u=0,\qquad a(t+T)=a(t)
\end{equation}
with $a\in C^1(\R)$, $0<a_1\le a(t)\le a_2$ and $\partial_t a \not\equiv 0$. Regardless of the `size' of the coefficient, this destroys any reasonable energy estimate as pointed out by Yagdjian, \cite{Yagdjian:2001}. We will give a short outline of the argument. Applying a partial Fourier transform to \eqref{eq:CPbeta} yields Hill's equation
\begin{equation}\label{eq:Hill}
  \hat u_{tt} + \lambda a^2(t) \hat u = 0
\end{equation}
with spectral parameter $\lambda=|\xi|^2\in [0,\infty)$. Due to a variant of Borg's theorem (see, e.g., the text book of Magnus-Winkler, \cite{Magnus:1966} or Colombini-Spagnolo \cite{Colombini:1984} for a refinement) there exists an open interval $\mathcal I\subset [0,\infty)$ such that for all $\lambda\in\mathcal I$ equation \eqref{eq:Hill} has an exponentially increasing solution. Taking appropriate Cauchy data for \eqref{eq:CPbeta} with Fourier support on $\mathcal I$ we can construct solutions with exponentially growing energy.

\begin{thm}[Yagdjian, \cite{Yagdjian:2001}] Assume $a\in L^1_{loc}(\R)$ is non-constant, periodic and a.e.{}~positive. Then
there exist Cauchy data $u(0,\cdot)=u_1$ and $u_t(0,\cdot)=u_2$ from Schwartz class such that the energy of the solution to \eqref{eq:CPbeta} satisfies
\begin{equation}
   \liminf_{t\to\infty} \frac{\log \E(u;t)}{\log t}  = \infty.
\end{equation}
\end{thm}

One can say a bit more about the structure of solutions to \eqref{eq:Hill}. There exist nested sequences $\lambda_j^\pm$, $0=\lambda_0^+<\lambda_1^-\le \lambda_1^+< \lambda_2^-\le \lambda_2^+<\cdots$, such that for all $\lambda\in (\lambda_j^+,\lambda_{j+1}^-)$ solutions are quasi-periodic, while for $\lambda\in\mathcal I_j=(\lambda_j^-,\lambda_j^+)$ solutions are combinations of exponentially increasing and exponentially decreasing contributions. Some of the intervals $\mathcal I_j$ may be empty, but not all.

\medskip
\noindent{\bf I.($\gamma$)} We can not ask for conservation of energy if we consider variable propagation speed (and no compensating lower order terms). However, we may ask for conditions on the coefficient function $a(t)$ such that the solutions to
\begin{equation}\label{eq:CPgamma}
  u_{tt} - a^2(t) \Delta u=0,\qquad u(0,\cdot)=u_1,\quad u_t(0,\cdot)=u_2
\end{equation}
satisfy Strichartz type decay estimates (and therefore also uniform energy bounds). A first answer was given by Reissig-Smith in \cite{Reissig:2005}. We assume
\begin{subequations}\label{eq:RS-cond}
\begin{align}
 & 0<a_1\le a(t)\le a_2,\\
 & \left| \frac{\d^k}{\d t^k} a(t) \right| \le C_k \left(\frac1{1+t}\right)^{k},\qquad k=1,2,\ldots,n+1.  \label{eq:gamma-symbest}
\end{align}
\end{subequations}
Then the following statement holds true:

\begin{thm}[Reissig-Smith, \cite{Reissig:2005}]\label{thm:2}
Assume $a\in C^{n+1}(\R)$ satisfies \eqref{eq:RS-cond}. Then the solution to \eqref{eq:CPgamma} satisfies the a priori estimate
\begin{equation}\label{eq:RS-strichartz}
  \| \nabla u(t,\cdot) \|_{L^q} + \| u_t(t,\cdot)\|_{L^q} \le C'_{pq} (1+t)^{-\frac{n-1}2 (\frac1p-\frac1q)} \big( \| \langle\D\rangle^{r_p+1} u_1\|_{L^p} + \| \langle\D\rangle^{r_p} u_2\|_{L^p} \big)
\end{equation}
for $p\in(1,2]$, $pq=p+q$ and $r_p = n(1/p-1/q)$.
\end{thm}

If we are interested in the energy alone, the constructed representation of solutions implies the two-sided energy estimate $C_1\E(u;0) \le \E (u;t) \le C_2 \E(u;0)$ for certain constants $C_1$ and $C_2$. For this result, \eqref{eq:gamma-symbest}  suffices for $k=1,2$. 

Condition \eqref{eq:gamma-symbest} can be weakened by the introduction of further $\log$-terms (but increasing also  the necessary number of derivatives). Then the estimate \eqref{eq:RS-strichartz} loses a small, resp. finite, amount of decay. For details we refer to \cite{Reissig:2005}. Results are sharp, there exist counter-examples to the estimate \eqref{eq:RS-strichartz} if the conditions are slightly violated.

\begin{expl}
  We give an example for an admissible coefficient function. The theory of Reissig-Smith allows to consider 
  $a(t)=2+\sin(\log t)$. For this coefficient the solutions to \eqref{eq:CPgamma} satisfy Strichartz-type decay estimates.
\end{expl}

The proof of Theorem~\ref{thm:2} is based on an explicit construction of the representation of solutions in terms of the coefficient function $a(t)$. Properties of the representation are different in different parts of the phase space $\R_t\times\R^n_\xi$, within the hyperbolic zone $\{(1+t)|\xi|\gg1\}$ it can be shown that solutions are represented by Fourier integrals
\begin{equation}\label{eq:FInt}
 {\sqrt{a(t)}}  \int \mathrm e^{\mathrm i (x\xi \pm |\xi| \int_0^t a(\theta)\d\theta)}  A(t,\xi)  \hat f(\xi)\d\xi, 
\end{equation}
$f$ given in terms of the initial data and $A$ being a symbol of order zero uniform in $t$. The construction is based on a diagonalisation argument within the hyperbolic zone. 

\medskip
\noindent{\bf I.($\delta$)} The `gap' between ($\beta$) and ($\gamma$) can be partially filled. Hirosawa, \cite{Hirosawa:2007}, weakened the symbol-like condition \eqref{eq:gamma-symbest} to smaller improvement per derivative in combination with a new assumption: the stabilisation condition. We consider the Cauchy problem \eqref{eq:CPgamma} and assume for the coefficient
\begin{subequations}\label{eq:H-cond}
\begin{align}
 & 0<a_1\le a(t)\le a_2,\\
 & \int_0^t |a(s)-a_\infty|\d s\le C (1+t)^q,\\
 & \left| \frac{\d^k}{\d t^k} a(t) \right| \le C_k \left(\frac1{1+t}\right)^{kp},\qquad k=1,2,\ldots, m.  \label{eq:delta-symbest}
\end{align}
\end{subequations}
with $m\ge 2$, $q\in[0,1)$ and $p\ge  q+(1-q)/m$. Then the following statement holds true:

\begin{thm}[Hirosawa, \cite{Hirosawa:2007}]\label{thm:3}
Assume $a\in C^m(\R)$ satisfies \eqref{eq:H-cond}. Then there exist constants $C_1,C_2>0$ such that the
generalised energy inequality
\begin{equation}
  C_1\E(u;0) \le \E (u;t) \le C_2 \E(u;0)
\end{equation}
is valid.
\end{thm}

On the contrary, there exists a coefficient satisfying \eqref{eq:H-cond} with $p<q$ for which no uniform energy bound holds true. Therefore, also this result is sharp. 

\begin{expl}
$a(t)=2+\sin(t^\alpha)$, $\alpha>0$, does not satisfy \eqref{eq:H-cond}. The stabilisation condition is violated. Indeed, with $a_\infty=2$ we obtain $\int_0^t |\sin(s^\alpha)|\d s=\frac1\alpha \int_0^{t^\alpha} |\sin(\theta)| \theta^{1/\alpha-1}\d \theta \simeq t$. 
\end{expl}

\begin{expl}
To obtain an admissible coefficient we follow \cite{Hirosawa:2007a} and choose sequences $\eta_j$, $\delta_j$ and $t_j$, $t_j+\delta_j<t_{j+1}$, together with a function $\psi\in C_0^\infty(\R)$ subject to $\mathrm{supp}\, \psi\subseteq[0,1]$, $-1<\psi(t)< 1$ and $\int_0^1|\psi(t)|\d t=1/2$. Then
\begin{equation}
  a(t) = 1 + \sum_{j=1}^\infty \eta_j \psi\left(\frac{t-t_j}{\delta_j}\right)
\end{equation}
satisfies \eqref{eq:H-cond} provided that $\eta_j\le 1$, $\sum_{j=1}^k \eta_j\delta_j \le C t_{j+1}^q$ and $t_j^p\le C \delta_j$. Such sequences exist, e.g., $t_j=2^j$, $\delta_j=2^{jq}$ and $\eta_j=2^{j(q-p)}$.
\end{expl}

The proof of Theorem~\ref{thm:3} is again based on explicit constructions of solutions. There are two major differences to the situation of Thm.~\ref{thm:2}. The hyperbolic zone is smaller $\{(1+t)^q|\xi|\gg 1\}$ (which makes it necessary to invoke a new argument for small frequencies) and the weaker assumptions for derivatives make it necessary to perform more diagonalisation steps. In consequence, representations by Fourier integrals contain an inhomogeneous phase function.

\begin{probl} Assume \eqref{eq:H-cond} for $p>q$ and all $k$. Is is possible to derive Strichartz type decay estimates under these assumptions? \end{probl}

\bigskip
{\bf II. Equations with increasing speed.} Problems with increasing speed have been considered by Reissig-Yagdjian, \cite{Reissig:2000}. Following their approach we write the coefficient function $a(t)$ as 
\begin{equation}
a(t) = \lambda(t) \omega(t)
\end{equation}
with a monotonously increasing function $\lambda(t)$ and a (bounded) oscillating part $\omega(t)$. We denote a primitive of $\lambda(t)$ as $\Lambda(t)=1+\int_0^t \lambda(s)\d s$ and assume a two-sided estimate
\begin{equation}
 \lambda'(t)\Lambda(t)\simeq(\lambda(t))^2.
\end{equation}
We refer to such functions $\lambda(t)$ as admissible shape functions and investigate the Cauchy problem
\begin{equation} \label{eq:CP-II}
u_{tt}-\lambda^2(t)\omega^2(t)\Delta u =0,\qquad u(0,\cdot)=u_1,\quad u_t(0,\cdot)=u_2
\end{equation}
for suitable $\omega(t)$.

\medskip
{\bf II.($\alpha$)} We shortly recall the result of Reissig-Yagdjian, \cite{Reissig:2000}, again neglecting logarithmic terms and the related classification of oscillations to simplify statements. We fix an admissible shape function $\lambda(t)$ and assume
\begin{subequations}\label{eq:RY-cond}
\begin{align}
 & 0<c_1 \lambda(t)\le a(t)\le c_2\lambda(t)\\
 &  \left| \frac{\d^k}{\d t^k} a(t) \right| \le C_k \lambda(t) \left(\frac{\lambda(t)}{\Lambda(t)}\right)^{k},\qquad k=1,2,\ldots, n+1, \label{eq:IIa-symbest}
\end{align}
\end{subequations}
Then the following Strichartz type decay estimate holds true:
\begin{thm}[Reissig-Yagdjian, \cite{Reissig:2000}]
Assume $a\in C^{n+1}(\R)$ satisfies \eqref{eq:RY-cond}. Then the solution to \eqref{eq:CP-II} satisfies the a priori estimate
\begin{equation}
  \| \lambda(t) \nabla u(t,\cdot) \|_{L^q} + \| u_t(t,\cdot)\|_{L^q} \le C'_{pq}{\sqrt{\lambda(t)}}\,(\Lambda(t))^{-\frac{n-1}2 (\frac1p-\frac1q)} \big( \| \langle\D\rangle^{r_p+1} u_1\|_{L^p} + \| \langle\D\rangle^{r_p} u_2\|_{L^p} \big)
\end{equation}
for $p\in(1,2]$, $pq=p+q$ and $r_p = n(1/p-1/q)$.
\end{thm}
If we consider the energy estimate alone, this yields $\E_\lambda(u;t) \le \lambda(t) \big(\|u_1\|_{H^1}^2+\|u_2\|^2_{L^2}\big)$ for the adapted energy
\begin{equation}
  \E_\lambda(u;t) = \frac12 \int \big( |\lambda(t) \nabla u(t,x)|^2+|u_t(t,x)|^2\big) \d x.
\end{equation} 

\medskip
{\bf II.($\beta$)} Again, we want to weaken the conditions to allow stronger oscillations in combination with a stabilisation condition in the spirit of \cite{Hirosawa:2007}. We assume that $a(t)=\lambda(t)\omega(t)$ satisfies
\begin{subequations}\label{eq:HW-cond}
\begin{align}
  &\lambda'(t) \simeq  \lambda(t) \left(\frac{\lambda(t)}{\Lambda(t)}\right),\quad
  |\lambda''(t)| \lesssim \lambda(t) \left(\frac{\lambda(t)}{\Lambda(t)}\right)^2, \quad \limsup_{t\to\infty} \frac{\lambda'(t)\Lambda(t)}{\lambda^2(t)} <2,\\
  & 0 < c_1\le \omega(t)\le c_2,
  \intertext{the stabilisation condition}
  & \int_0^t \lambda(s) |\omega(s)-\omega_\infty| \d s\le C \Theta(t) = o(\Lambda(t))
  \intertext{for an auxiliary function $\Theta(t)\le \Lambda(t)$, $\Theta(0)=1$, and a constant $\omega_\infty$,  the symbol-like estimate} 
  &  \left| \frac{\d^k}{\d t^k} a(t) \right| \le C_k \lambda(t) \left(\Xi(t)\right)^{-k},\qquad k=1,2,\ldots, m 
  \intertext{with an (increasing) function $\Xi(t)$, $\lambda(t)\Xi(t)\ge C\Theta(t)$, and}
  & \int_t^\infty (\lambda(s))^{1-m} (\Xi(s))^{-m}\d s\le C (\Theta(t))^{1-m}.
\end{align}
\end{subequations}
Under these assumptions the following energy estimate holds true:

\begin{thm}[Hirosawa-Wirth, \cite{Hirosawa:2007a}]
Assume $a\in C^m(\R)$ satisfies \eqref{eq:HW-cond}. Then solutions to the Cauchy problem \eqref{eq:CP-II} satisfy the energy estimate
\begin{equation}
\E_\lambda(u;t) \le \lambda(t) \big(\|u_1\|_{H^1}+\|u_2\|_{L^2}\big).
\end{equation} 
Furthermore, if $(u_1,u_2)\not\equiv (0,0)$ the limit
\begin{equation}
  \lim_{t\to\infty} \frac1{\lambda(t)}  \E_\lambda(u;t) \ne 0
\end{equation}
exists and is non-zero.
\end{thm}

\begin{probl}
Derive dispersive estimates under similar assumptions. This should be a consequence of the solution to Problem~1.
\end{probl}

\bigskip
{\bf III. Equations with non-effective dissipation.} We consider the influence of dissipation terms and investigate
\begin{equation}\label{eq:CP-III}
  u_{tt}-\Delta u+2b(t)u_t = 0,\qquad u(0,\cdot)=u_1,\quad u_t(0,\cdot)=u_2.
\end{equation}
This Cauchy problem is related to \eqref{eq:CP-II} with increasing propagation speed by a change of variables. For completeness we give the correspondence. Let $v(t,x)$ be a solution to $v_{tt}-\lambda^2(t)\Delta v = 0$. Then 
$u(t,x) = v(\Lambda^{-1}(t),x)$---with $\Lambda^{-1}(t)$ the inverse function to $\Lambda(t)$---satisfies
\begin{subequations}
\begin{align}
  &u_t = \frac{v_t(\Lambda^{-1}(t),x)}{ \lambda(\Lambda^{-1}(t))},\\
  &u_{tt} = \frac{v_{tt}(\Lambda^{-1}(t),x)}{\lambda^2(\Lambda^{-1}(t))}-\frac{\lambda'(\Lambda^{-1}(t))}{\lambda^3(\Lambda^{-1}(t))}  v_t(\Lambda^{-1}(t),x) = 
  \Delta u -   \frac{\lambda'(\Lambda^{-1}(t))}{\lambda^2(\Lambda^{-1}(t))} u_t , 
\end{align}
\end{subequations}
which is of the form \eqref{eq:CP-III}; $\lambda(t)= (1+t)^\ell$ corresponds asymptotically to $2b(t)\sim \frac{\ell}{\ell+1} \frac1{1+t}$.

\medskip
{\bf III.($\alpha$)} In \cite{Wirth:2004c}, \cite{Wirth:2004}, \cite{Wirth:2006} the author gave an overview on energy and dispersive estimates for solutions to \eqref{eq:CP-III} under assumptions related to \cite{Reissig:2000}, \cite{Reissig:2005}. We will quote a corresponding result. We assume for the coefficient function
\begin{subequations}\label{eq:W-cond}
\begin{align}
   & b(t)\ge 0, \\
   & \limsup_{t\to\infty} tb(t) < \frac12,\\ 
   & \left|\frac{\d^k}{\d t^k} b(t)\right| \le C_k \left(\frac1{1+t}\right)^{1+k},\qquad k=1,2,\ldots, n+2.
\end{align}
\end{subequations}
Then the following Strichartz type decay estimate holds true:
\begin{thm}[Wirth, \cite{Wirth:2006}]
Assume $b\in C^{n+2}(\R)$ satisfies \eqref{eq:W-cond}. Then the solution to \eqref{eq:CP-III} satisfies the a priori estimate
\begin{equation}\label{eq:W-strichartz}
  \| \nabla u(t,\cdot) \|_{L^q} + \| u_t(t,\cdot)\|_{L^q} \le C''_{pq} \frac1{\beta(t)} (1+t)^{-\frac{n-1}2 (\frac1p-\frac1q)} \big( \| \langle\D\rangle^{r_p+1} u_1\|_{L^p} + \| \langle\D\rangle^{r_p} u_2\|_{L^p} \big)
\end{equation}
for $p\in(1,2]$, $pq=p+q$ and $r_p = n(1/p-1/q)$ and with the auxiliary function
\begin{equation}
  \beta(t) = \exp\int_0^t b(s)\d s.
\end{equation}
Furthermore, for non-zero initial data $(u_1, u_2)\in H^1\times L^2$ the limit
\begin{equation}\label{eq:scatt-est}
  \lim_{t\to\infty} \beta^2(t)\E(u;t) \ne 0
\end{equation}
exists and is non-zero.
\end{thm}

In particular, \eqref{eq:scatt-est} can only be true, if $tb(t)$ is asymptotically small. As pointed out by
Hirosawa-Nakazawa, \cite{Hirosawa:2003a}, the example $b(t)=\mu / (1+t)$ with $\mu>1$ satisfies
\begin{equation}\label{eq:HN-est}
  \lim_{t\to\infty} t^2 \E(u;t) = 0
\end{equation}
for all data and, furthermore, the exponent $2$ is best possible (see, e.g., \cite{Wirth:2004c}). 

\medskip
{\bf III.($\beta$)} The conditions \eqref{eq:W-cond} can be weakened in the spirit of \cite{Hirosawa:2007} to allow faster oscillations of the coefficient. We assume that $2b(t) = \mu(t) + \sigma(t)$, where $\mu(t)$ is a monotonically decreasing shape function and $\sigma(t)$ carries the oscillations. Conditions are similar to II.($\beta$) and follow
\cite{Hirosawa:2007b}. We assume
\begin{subequations}\label{eq:HW2-cond}
\begin{align}
 & \mu(t)>0, \quad \mu'(t)<0,\quad \limsup_{t\to\infty} t\mu(t) < 1, \\
 & \sup_{t>0} \left| \int_0^t \sigma(s)\d s\right| < \infty,\\
 & \int_0^t \left| \exp\left( \int_0^\theta \sigma(s)\d s\right) - \omega_\infty \right| \d\theta \le C \Theta(t),
 \intertext{with an auxiliary function $\Theta(t)=o(t)$, $t\to\infty$,}
 & \left| \frac{\d^k}{\d t^k} b(t) \right| \le C_k (\Xi(t))^{-k-1}, \qquad k=0,1,2,\ldots, m
 \intertext{with $\Xi(t) \gtrsim \Theta(t)$ and}
 &\int_t^\infty (\Xi(s))^{-m-1}\d s\lesssim (\Theta(t))^{-m}.
\end{align}
\end{subequations}
 Then the following statement holds true:
 
 \begin{thm}[Hirosawa-Wirth,\cite{Hirosawa:2007b}]
 Assume $b\in C^{m}(\R)$ satisfies \eqref{eq:HW2-cond}. Then the energy of the solution to \eqref{eq:CP-III} satisfies
 \begin{equation}
  \E(u;t) \le C \frac1{\beta^2(t)} \left( \|u_1\|_{H^1} + \|u_2\|_{L^2} \right)
 \end{equation}
 with $\beta(t)=\exp\int_0^t b(s)\d s\simeq \exp\int_0^t \mu(s)\d s$. Furthermore, for non-zero initial data the limit
 \begin{equation}
   \lim_{t\to\infty} \beta^2(t) \E(u;t) \ne 0
 \end{equation}
 exists and is non-zero.
 \end{thm}

\begin{expl}
We can use $\mu(t) = \frac\mu{1+t}$ and $\sigma(t) = \mu(t) \sin(t^\alpha)$ with $\alpha\in(0,1)$. Then all conditions are satisfied. Note, that the corresponding approach due to Reissig-Smith allows only $\sin(\log t)$ terms.
\end{expl}

\begin{expl}
We can come closer to $\sin(t)$ by using $\mu(t)=\frac1{(1+t)\log(e+t)}$ and $\sigma(t) = \mu(t) \sin(t/\log(e+t))$.
\end{expl}

Note, that $\sigma(t)$ has \emph{no} influence on the decay rates, only on appearing constants. $\sigma(t)$ can not be seen as a small perturbation of $\mu(t)$; it might be arbitrary large. 

\bigskip
{\bf IV. Equations with effective dissipation.} Following the classification of \cite{Wirth:2004}, \cite{Wirth:2007}, we speak of effective dissipation if the lower order term $2b(t)u_t$ changes the large time asymptotics of solutions to \eqref{eq:CP-III} in an essential way. 

\medskip
{\bf IV.($\alpha$)} The damped wave or telegraph equation 
\begin{equation}\label{eq:CP-damped}
  u_{tt} - \Delta u + u_t =0,\qquad u(0,\cdot)=u_1, \quad u_t(0,\cdot)=u_2
\end{equation}
was studied by Matsumura in \cite{Matsumura:1976}. He investigated semilinear perturbations of this equation based on
\begin{thm}[Matsumura, \cite{Matsumura:1976}]
  The solutions to the damped wave equation \eqref{eq:CP-damped} satisfy the a-priori estimates
  \begin{subequations}
  \begin{align}
   & \|\partial_t^k \partial_x^\alpha u(t,x) \|_{L^2} \le C (1+t)^{-\frac{|\alpha|}2-k-\frac n4} \big( \|u_1\|_{H^{k+|\alpha|}} + \|u_2\|_{H^{k+|\alpha|-1}} \notag\\
   &\qquad\qquad\qquad\qquad\qquad\qquad\qquad\qquad\qquad\qquad+ \|u_1\|_{L^1} + \|u_2\|_{L^1}  \big) \\
    & \|\partial_t^k \partial_x^\alpha u(t,x) \|_{L^\infty} \le C (1+t)^{-\frac{|\alpha|}2-k-\frac n2} \big( \|u_1\|_{H^{r+k+|\alpha|}} + \|u_2\|_{H^{r+k+|\alpha|-1}}  \notag\\
   &\qquad\qquad\qquad\qquad\qquad\qquad\qquad\qquad\qquad\qquad+ \|u_1\|_{L^1} + \|u_2\|_{L^1}  \big)
  \end{align}
  \end{subequations}
  with $r>\frac n2$.
\end{thm}

The estimates are different in structure compared to the Strichartz type estimates for the free wave equation and the related estimates discussed so far. While the used regularity is related to Sobolev embedding (and therefore hyperbolic in nature) the decay rates correspond to estimates for the heat equation (as, e.g., used by Ponce in \cite{Ponce:1985}). This parabolic structure becomes even more apparent in the diffusion phenomenon due to
Nishihara, \cite{Nishihara:1997}, \cite{Nishihara:2003}, and Yang-Milani, \cite{Yang:2000} or the estimates of Narazaki, \cite{Narazaki:2004}.

\begin{thm}[Nishihara, \cite{Nishihara:2003}] 
Let $u(t,x)$ be a solution to \eqref{eq:CP-damped} in $\R_t\times\R^3_x$, $w(t,x)$ the corresponding solution to
\begin{equation}
  w_t=\Delta w,\qquad w(0,\cdot)=u_1+u_2
\end{equation}
and $v(t,x)$ a (suitable) free wave. Then the $L^p$--$L^q$ estimate
\begin{equation}
  \|u(t,\cdot)-w(t,\cdot)-\mathrm e^{-\frac t2} v(t,\cdot)\|_{L^q}\le  C (1+t)^{-\frac32(\frac1p-\frac1q)-1} \left(\|u_1\|_{L^p} + \|u_2\|_{L^p}\right)
\end{equation}
holds true for all $1\le p\le q\le\infty$.
\end{thm}

The decay rate is by one order better than expected for both $\|u(t,\cdot)\|_{L^q}$ and $\|w(t,\cdot)\|_{L^q}$ (while the regularity of the data on the right hand side is not sufficient to obtain these decay estimates---singularities are eliminated by the free wave $v(t,x)$). 

For both statements it is essential that {\em small} frequencies determine the asymptotic behaviour of solutions. If we localise to high frequencies (or to a hyperbolic part of the phase space), we obtain exponential decay rather than the polynomial decay given in the theorems. 

\medskip
{\bf IV.($\beta$)} In \cite{Wirth:2007} the author discussed variable coefficient versions of these results; we quote a slightly simplified version (c.f. \cite{Reissig:2006} or \cite{Wirth:2004} for this formulation). Let $b(t)$ be subject to
\begin{subequations}\label{eq:RW-cond}
\begin{align}
&b(t)\ge 0,\qquad b'(t)\ne 0,\qquad \lim_{t\to\infty} tb(t)=\infty,\\
&\left|\frac{\d^k}{\d t^k} b(t)\right| \le C_k b(t) \left(\frac1{1+t}\right)^k,\qquad k=1,2.
\end{align}
\end{subequations}
Then the following decay estimate holds true:
\begin{thm}[Wirth, \cite{Wirth:2004}] 
  Assume $b\in C^{2}(\R)$ satisfies \eqref{eq:RW-cond}. Then the solution to \eqref{eq:CP-III} satisfies the a priori estimate
  \begin{subequations}\label{eq:RW-strichartz}
\begin{align}
  \| u(t,\cdot) \|_{L^q}  &\le C_{pqr} \left(1+\int_0^t \frac{\d s}{b(s)}\right)^{-\frac{n}2 (\frac1p-\frac1q)} \big( \| \langle\D\rangle^{r_p} u_1\|_{L^p} + \| \langle\D\rangle^{r_p-1} u_2\|_{L^p} \big)\\
  \| \nabla u(t,\cdot) \|_{L^q}  &\le C_{pqr} \left(1+\int_0^t \frac{\d s}{b(s)}\right)^{-\frac{n}2 (\frac1p-\frac1q)-\frac12} \big( \| \langle\D\rangle^{r_p+1} u_1\|_{L^p} + \| \langle\D\rangle^{r_p} u_2\|_{L^p} \big)\\
  \| u_t(t,\cdot)\|_{L^q}&\le C_{pqr} \left(1+\int_0^t \frac{\d s}{b(s)}\right)^{-\frac{n}2 (\frac1p-\frac1q)-1} \big( \| \langle\D\rangle^{r_p+1} u_1\|_{L^p} + \| \langle\D\rangle^{r_p} u_2\|_{L^p} \big)
\end{align}
\end{subequations}
for $1\le p\le 2\le q\le\infty$ and $r_p > n(1/p-1/q)$.
\end{thm}

The proof of this statement is constructive in the sense that the main terms in the representation of solutions are given explicitly. The assumption $tb(t)\to\infty$ guarantees that the high frequency part decays faster than the low frequency part and is essential for the parabolic type estimates. Under slightly more restrictive assumptions the diffusion phenomenon can be obtained. We refer to \cite{Wirth:2007} for the details.

\medskip
{\bf IV.($\gamma$)} We conclude the overview with some remarks concerning periodic dissipation terms. In contrast to periodic propagation speed this does not destroy Matsumura's decay estimates nor the diffusion phenomenon. The following two statements are based on \cite{Wirth:2007d}.

\begin{thm}[Wirth, \cite{Wirth:2007d}]
Assume $b\in AC_{loc}(\R)$ is absolutely continuous, positive a.e.~and periodic. Then the solution to \eqref{eq:CP-III} satisfies the a priori estimate
\begin{subequations}\label{eq:W2-strichartz}
\begin{align}
 \| u(t,\cdot) \|_{L^q} &\le C'_{pqr} (1+t)^{-\frac{n}2 (\frac1p-\frac1q)} \big( \| \langle\D\rangle^{r_p} u_1\|_{L^p} + \| \langle\D\rangle^{r_p-1} u_2\|_{L^p} \big)\\
  \| \nabla u(t,\cdot) \|_{L^q} &\le C'_{pqr} (1+t)^{-\frac{n}2 (\frac1p-\frac1q)-\frac12} \big( \| \langle\D\rangle^{r_p+1} u_1\|_{L^p} + \| \langle\D\rangle^{r_p} u_2\|_{L^p} \big)\\
   \| u_t(t,\cdot)\|_{L^q}&\le C'_{pqr} (1+t)^{-\frac{n}2 (\frac1p-\frac1q)-1} \big( \| \langle\D\rangle^{r_p+1} u_1\|_{L^p}  \| \langle\D\rangle^{r_p} u_2\|_{L^p} \big)
  \end{align}
\end{subequations}
for $1\le p\le2\le q\le\infty$ and $r_p = n(1/p-1/q)$.
\end{thm}

\begin{probl} Close the gap between IV.($\beta$) and IV.($\gamma$), i.e., develop energy estimates based on less regularity (and/or less improvement per derivative).  
\end{probl}

\begin{thm}[Wirth, \cite{Wirth:2007d}]
 Assume $b\in AC_{loc}(\R)$ is absolutely continuous, positive a.e.~and periodic. Then there exist constants $\alpha$ and $\beta$ depending only on the coefficient function $b$ such that the solution to \eqref{eq:CP-III} and the corresponding solution to 
 \begin{equation}
    w_t=\alpha\Delta w,\qquad w(0,\cdot)=u_1+\beta u_2
 \end{equation}
 satisfies
 \begin{equation}
   \|u(t,\cdot)-w(t,\cdot)\|_{L^2} \le C (1+t)^{-1} \big(\|u_1\|_{H^1} + \|u_2\|_{L^2}\big).
 \end{equation}
\end{thm}

The constants $\alpha$ and $\beta$ are explicit and given in \cite{Wirth:2007d}. A similar estimate holds true if we replace the $L^2$-norm on the right hand side by an $L^q$-norm, $q\ge 2$, or measure spatial derivatives of the difference.

Both theorems are essentially based on the idea that small frequencies determine the asymptotic behaviour. Although the periodic coefficient leads to a variant of Hill's equation after partial Fourier transform and a change of variables, we only need to understand the neighbourhood of $\xi=0$. It can be shown that away from this exceptional frequency, solutions decay exponentially, while nearby (which is inside an instability interval) the structure of the fundamental solution is known explicitly due to Floquet's theory (cf., \cite{Magnus:1966} or \cite[Chapter IXX]{Watson:1922}). 

\bigskip
{\bf V. Concluding remarks.} We focussed on the influence of variable speed and variable dissipation. There are further results on the treatment of variable mass terms, see e.g. the treatment of B\"ohme-Reissig, \cite{Bohme:2009}, of scale invariant mass or \cite{Hirosawa:2003} on interactions between variable mass and variable speed. 

There is a fundamental difference between oscillations in the principle part and oscillations in lower order terms; mainly due to asymptotic relations to problems of different type. While in the first case sharp results are available (within certain symbol-like classes of coefficients), the example of effective dissipation shows an essential gap of understanding between \cite{Wirth:2007} and \cite{Wirth:2007d}, which needs to be closed in the future.  


\end{document}